\documentclass{amsart}

\newtheorem{theorem}{Theorem}[section]
\newtheorem{proposition}[theorem]{Proposition}

\newtheorem{definition}[theorem]{Definition}

\def\a{\alpha}
\def\b{\beta}
\def\d{\delta}

\def\F{\Phi}
\def\g{\gamma}

\def\k{\kappa}
\def\l{\lambda}
\def\L{\Lambda}
\def\o{\omega}
\def\O{\Omega}

\def\s{\sigma}
\def\S{\Sigma}
\def\t{\tau}
\def\te{\theta}

\def\z{\zeta}

\def\ov{\overline}

\def\wt{\widetilde}
\def\wh{\widehat}

\def\SZ{Szeg\H o}
\def\({\left(}
\def\){\right)}
\def\[{\left[}
\def\]{\right]}

      \def\dC{{\mathbb C}}
\def\dD{{\mathbb D}}

   \def\dN{{\mathbb N}}   \def\dO{{\mathbb O}}
\def\dP{{\mathbb P}}      \def\dR{{\mathbb R}}

   \def\dZ{{\mathbb Z}}

      \def\cC{{\mathcal C}}

      \def\cL{{\mathcal L}}
\def\cM{{\mathcal M}}

\def\ds{\displaystyle}

\font\tenopen = cmbx10 \font\sevenopen = cmbx7 \font\fiveopen =
cmbx5
\newfam\openfam

\textfont\openfam = \tenopen \scriptfont\openfam = \sevenopen
\scriptscriptfont\openfam = \fiveopen

\begin{document}
\title[Rational interpolation and mixed inverse problem]
{Rational interpolation and mixed inverse spectral problem for
finite CMV matrices}

\author{Leonid Golinskii}
\address{Institute for Low Temperature Physics and Engineering\\47\\ Lenin
ave.\\ Kharkov \\ Ukraine} \email{golinsky@ilt.kharkov.ua}

\author{Mikhail Kudryavtsev}
\address{Institute for Low Temperature Physics and Engineering\\47\\ Lenin
ave.\\ Kharkov \\ Ukraine} \email{kudryavtsev@ilt.kharkov.ua;
kudryavstev@onet.com.ua}



\begin{abstract}
For finite dimensional CMV matrices the mixed inverse spectral
problem of reconstruction the matrix by its submatrix and a part
of its spectrum is considered. A general rational interpolation
problem which arises in solving the mixed inverse spectral problem
is studied, and the description of the space of its solutions is
given. We apply the developed technique to give sufficient
conditions for the uniqueness of the solution of the mixed inverse
spectral problem.
\end{abstract}

\maketitle

 \section{Introduction}\setcounter{section}{1}
\setcounter{equation}{0}

The theory of CMV matrices, rapidly developing in the recent years
(\cite{CMV2, CMV3, CMV4}, see also the expositions in \cite{sim2,
KN} and the references therein), has a strong background in the
theory of Jacobi matrices. The similarity between the Jacobi and
CMV matrices not only provides a general concern of investigation,
but also permits sometimes to predict the most probable answers.
Some of the facts known for the Jacobi matrices can be easily
extended to the CMV case, some require considerable adjustment. In
turn, the spectral theory of the Jacobi matrices is paralleled by
the spectral theory for the  Sturm-Liouville differential
operators, which forms another important front of study in this
area of research.

{\it Mixed inverse spectral problems} ({\bf MISP}) are of special
interest in the inverse spectral theory. For this kind of problems
one reconstructs a differential or difference operator by a part of
its potential and some additional spectral data. Started by
Hochstadt and Lieberman \cite{HL} for the Sturm-Liouville operators,
these problems were extended and refined for many other cases (see
\cite{Chu-Golub} for the references). Compared to the ``ordinary''
inverse spectral problems where the whole potential is to be
reconstructed, for MISP one needs to know ``less'' spectral data.

In what follows, we refer to several key studies of the
development of MISP. Let $J$ be an $N\times N$ three-diagonal
matrix of the form
$$
J=\begin{pmatrix} a_1 & b_1 & 0 & 0 & \cdot & \cdot & \cdot \\
b_1 & a_2 & b_2 & 0 & \cdot & \cdot & \cdot \\
0 & b_2 & a_3 & b_3 & \cdot & \cdot & \cdot \\
\cdot & \cdot & \cdot & \cdot & \cdot & \cdot & \cdot \\
\cdot & \cdot & \cdot & \cdot & \cdot & \cdot & \cdot \\
\cdot & \cdot & \cdot & \cdot & 0 & b_{N-1} & a_N\\
\end{pmatrix}, \quad a_n, b_n \in {\dR}, \ b_n>0.
$$
Consider the $a$'s and the $b$'s as a single sequence $a_1, b_1,
a_2, b_2, \ldots$, that is,
$$ c_{2n-1}=a_n, \quad c_{2n}=b_n, \qquad n\in \dN.
$$
In \cite{H} Hochstadt proved the discrete version of the
Hochstadt-Lieberman theorem:

\medskip

{\bf Theorem.} {\it Let $n\in \dN$. Suppose that
$c_{N+1},\ldots,c_{2N-1}$ are known, as well as the eigenvalues
$\l_1,\ldots,\l_N$ of $J$. Then $c_1,\ldots,c_N$ are uniquely
determined.}

\medskip

It is important to remark that in certain complicated physical
systems it is not always possible to know the entire spectrum. A
natural question arises if it is possible to reconstruct the Jacobi
matrix when we know more than a half of the potential, but less than
the whole spectrum. The positive answer is given by Gesztesy and
Simon in \cite{GS}.

\medskip

{\bf Theorem.} {\it Suppose that $1\leq j\leq N$ and
$c_{j+1},\ldots,c_{2N-1}$ are known, as well as (any) $j$ of the
eigenvalues. Then $c_1,\ldots,c_j$ are uniquely determined.}

\medskip

The goal of this paper is to study the MISP for CMV matrices.
Although the algorithm of solving the problem is similar to that
of the Jacobi case, essential difference arises when the
uniqueness of the reconstruction is concerned. In the Jacobi case
the MISP is reduced to the interpolation of a rational function
(specifically, the Weyl function of the unknown submatrix) by its
values in the known eigenvalues of the whole matrix. The degree of
the rational function (i.e., the sum of the degrees of the
numerator and denominator) corresponds to the number of the
interpolation points, and both numerator and denominator are
monic. So, the uniqueness of the interpolating rational function
drops out immediately from the simple fact that a polynomial of
degree $k$ with $k+1$ zeros is identically zero, and thus we prove
the uniqueness of the reconstructed Jacobi matrix. However, in the
CMV case the degree of the interpolating rational function (the
Weyl function of the ``reduced CMV matrix'') is greater by $1$
than the number of the interpolation points, and the lacking
``piece of information'' is given by a restricting condition on
the free term of the numerator to be $1$. Here the trivial
consideration of the Jacobi case fails, and a theory is required
of how to find the interpolating rational function with such
restriction.

So, the work consists of two parts. In Section 2 we give an
approach to the rational interpolation theory adapted for solution
of the interpolation problem related to the MISP. The MISP itself
is studied in Section 3. Note that, although in the Jacobi case we
always have uniqueness in the MISP (provided the number of lacking
entries agrees with the number of the given parameters), in the
CMV case a degenerated case is possible where the MISP has
infinitely many solutions. The main results of this paper are
Theorems \ref{theor2.9} with a description of the solutions of the
rational interpolation problem, and \ref{theormain} with a
sufficient condition for the MISP to have a unique solution.

\setcounter{section}{1} \setcounter{equation}{0}
\section {Two-dimensional vector-polynomials and rational interpolation}

We start with the following interpolation problem: given points
$z_1,z_2,\ldots,z_n \in \mathbb{C}$ and numbers
$\omega_1,\omega_2,\ldots,\omega_n \in \overline{\mathbb{C}}$ find a
``nice'' description of {\it all} rational functions
$P^{(1)}/P^{(2)}$, where $P^{(1)}$ and $P^{(2)}$ are polynomials
with complex coefficients, for which
\begin{equation} \label{2.1}
{ P^{(1)}(z_j) \over P^{(2)}(z_j) } = \omega_j \,, \qquad
j=1,\ldots,n,
\end{equation}
($\omega_j=\infty$ means that the rational function must have a
pole at the point $z_j$). By a ``nice'' description we understand
a description whose form allow us to ``control'' the degrees of
$P^{(j)}(z)$, $j=1,2$, and some more parameters such as their
leading coefficients or the free term, required in the concrete
application od the rational interpolation problem (we will return
later on to this topic).

It is advisable to reformulate this problem as a linear problem in
the space of two-dimensional polynomial vector-functions
(vector-polynomials). Define numbers $\alpha^{(1)}_j$ and
$\alpha^{(2)}_j$ by the following rule

\begin{equation} \label{2.2}
\left\{\begin{array}{l} \alpha_{j,1} := 1, \quad \alpha_{j,2} :=
-\omega_j, \ \  {\rm if} \ \omega_j \ne \infty , \\
\alpha_{j,1} := 0, \quad \alpha_{j,2} := 1, \qquad {\rm if}
\ \omega_j = \infty .\\
\end{array}
\right.
\end{equation}
Then (\ref{2.1}) implies
\begin{equation} \label{2.3}
\alpha_{j,1} P^{(1)}(z_j) + \alpha_{j,2} P^{(2)}(z_j) = 0, \qquad
j=1,\ldots,n,
\end{equation}
where $|\alpha_{j,1}| + |\alpha_{j,2}| > 0$.

Conversely, if $|P^{(1)}(z_j)| + |P^{(2)}(z_j)| >0$, $j=1,\ldots,n$,
then (\ref{2.3}) implies (\ref{2.1}) with
\begin{equation} \label{2.4}
\omega_j := - \frac{\alpha_{j,1}}{\alpha_{j,2}} \in
\overline{\mathbb{C}}, \qquad j=1,\ldots,n.
\end{equation}
However, if $|P^{(1)}(z_j)| + |P^{(2)}(z_j)|=0$ for some $j$, then
(\ref{2.1}) does not make sense at the point $z_j$.

Thus, we reformulated the initial interpolation problem in the
following way: find all the pairs of polynomials $P^{(1)}$ and
$P^{(2)}$, for which (\ref{2.3}) holds. The second problem is
"almost equivalent" to the first one: if the rational function
$P^{(1)}/P^{(2)}$ is a solution of (\ref{2.1}), then the pair of
polynomials $P^{(1)}$ and $P^{(2)}$ is a solution of (\ref{2.3})
with $\alpha_{j,1}$ and $\alpha_{j,2}$ defined in (\ref{2.2}). If
$P^{(1)}$ and $P^{(2)}$ are polynomials such that $|P^{(1)}(z_j)| +
|P^{(2)}(z_j)|>0$ solving (\ref{2.3}), then the rational function
$P^{(1)}/P^{(2)}$ is a solution of (\ref{2.1}) with $\o_j$ defined
in (\ref{2.4}).

In what follows we refer to interpolation problem \eqref{2.3} as the
problem $(I_n)$, and denote by $\dP_\infty(I_n)$ the set of all its
solutions. Some of the results desribed below for the polynomial
vector-functions are taken from \cite{Kud1}, \cite{Kud2} where the
interpolation problems also appears as well as its main objects like
generators, etc.

\subsection{The space of vector-polynomials}

To describe the class of solutions of (\ref{2.3}), we introduce the
space
$$ \mathbb{P}_\infty := \left\{ p(z) =
\left(\begin{array}{c} P^{(1)}(z) \\ P^{(2)}(z)\\
\end{array}\right), \quad  P^{(i)}, \ i=1,2, \ {\rm
are \ complex \ polynomials}  \right\}.
$$
This is a linear space with the standard operations. The zero
element on
this space is $0 =\left(\begin{array}{c} 0 \\ 0 \\
\end{array}\right)$. We point out that $\dP_\infty$ is also a module over
the ring of polynomials:
$$ Sp=\left(\begin{array}{c} SP^{(1)} \\ SP^{(2)}\\
\end{array}\right)\in \dP_\infty $$
for any polynomial $S$.

\begin{definition}\label{def1}
The height of the vector-polynomial $p =
\left(\begin{array}{c} P^{(1)} \\ P^{(2)}\\
\end{array}\right) \ne 0$ is the number
\begin{equation} \label{2.5}
h(p) := \left\{ \begin{array}{l} 2 \deg P^{(1)}, \quad \deg P^{(1)}
> \deg P^{(2)},\\
2 \deg P^{(2)} + 1, \quad \deg P^{(1)}
\leq \deg P^{(2)},\\
\end{array}
\right.
\end{equation}
As usual, $\deg 0=-\infty$, so we put $h(0):= -\infty$.
\end{definition}

It is obvious that
\begin{equation} \label{2.6}
h(p) = \max \{ 2\deg P^{(1)}, 2\deg P^{(2)}+1 \}
\end{equation}
and
\begin{equation} \label{2.61}
h(Sp)=h(p)+2\deg S.
\end{equation}
 The degrees of the components of the vector-polynomials can be
written down in the following table:

\medskip

\begin{flushleft}
\begin{tabular}{|c|c|c|c|c|c|c|c|c|c}
\hline {height p} & $0$ & $1$ & $2$ & $3$ &
$4$ & \ldots & $2k$ & $2k+1$ & \ldots \\
\hline
$\deg P^{(1)}$ & $=0$ & $\leq 0$ & $=1$ & $\leq 1$ & $=2$ & \ldots & $=k$ & $\leq k$ & \ldots \\
\hline
$\deg P^{(2)}$ & $-\infty$ & $=0$ & $\leq 0$ & $=1$ & $\leq 1$ & \ldots & $\leq k-1$ & $=k$ & \ldots\\
\hline
\end{tabular}
\end{flushleft}

\medskip


The following proposition demonstrates that the notion of the height
is a natural extension of the degree of polynomials.

\begin{proposition} \label{propos2.1}
{\it If $h(p) \ne h(q)$, then $\forall a,b \in \mathbb{C}$
$$ h(ap + bq) = \max(h(p),h(q))
$$
 If $h(p) = h(q) = n$, then
 \begin{enumerate}
    \item    $\forall \ a,b \in \mathbb{C}$ \ \ $h(ap+bq) \leq n$
    \item    $\exists c \in \mathbb{C}$: \ \ \ \ \ $h(p+cq) \leq n-1$.
 \end{enumerate}
}\end{proposition}

\noindent {\it Proof.} We prove (2), the rest is plain. If $n=2k$,
then $k=\deg P^{(1)}>\deg P^{(2)}$, $k= \deg Q^{(1)}>\deg Q^{(2)}$.
So, there exists $c\in\mathbb{C}$ such that $\deg(P^{(1)} + c \,
Q^{(1)}) \leq k-1$. Also, since $\deg P^{(2)} \leq k-1$ and $\deg
Q^{(2)} \leq k-1$, we see that $\deg(P^{(2)} + c \, Q^{(2)}) \leq
k-1$. Then, by (\ref{2.6})
$$ h(p+cq) \leq 2k-1=n-1.
$$
If $n=2k+1$, then $k=\deg P^{(2)}\ge\deg P^{(1)}$, $k=\deg Q^{(2)}
\ge\deg Q^{(1)}$. So, $\exists c\in\mathbb{C}$: $\deg(P^{(2)}+c\,
Q^{(2)}) \leq k-1$, $\deg(P^{(1)} + c \, Q^{(1)}) \leq k$. Then, by
(\ref{2.6})
$$ h(p+cq) \leq 2k=n-1.
$$
\hfill $\square$

Consider the following basic system of vectors in
$\mathbb{P}_\infty$:
 \begin{equation} \label{2.62}
e_{2k}(z):=
\begin{pmatrix} z^k\\ 0\end{pmatrix}, \ \
e_{2k+1}(z):= \begin{pmatrix} 0\\ z^k
\end{pmatrix}, \quad k\in\dZ_+=\{0,1,\ldots\}
\end{equation}
It is clear that $h(e_n)=n$ for all $n$.

\begin{proposition} \label{propos2.2}
{\it $\{e_n\}_{n\ge0}$ is a basis in $\dP_\infty$, i.e.,
$$ p\in \dP_\infty, \quad h(p)=m \quad \Longrightarrow \quad
p(z)=\sum_{k=0}^mc_ke_k(z), \quad c_m\ne0,
$$
and this representation is unique.}
\end{proposition}

{\it Proof.} \ We use induction on $m$:

1) The cases $m=0,1$ are checked immediately

2) Let $m=2k+1$, then $k=\deg P^{(2)}\ge \deg P^{(1)}$, i.e.,
$P^{(2)}(z)=a z^k + Q^{(2)}(z)$, $a\ne0$, $\deg Q^{(2)}(z)\leq k-1$.
Then
$$ q(z):= p(z) - a e_{2k+1}(z)=
\begin{pmatrix} Q^{(1)}(z)\\ Q^{(2)}(z)
\end{pmatrix},
$$
where $Q^{(1)}(z)=P^{(1)}(z)$ and $\deg Q^{(1)}\leq k$.

So, $h(q)\leq 2k = m-1$, and we can apply the inductive hypothesis.
The representation $q(z)=\sum_{k=0}^{m-1}c_ke_k(z)$ is unique, so
the representation
$$ p(z)=a e_m(z) + \sum_{k=0}^{m-1}c_ke_k(z)
$$
is also unique.

3) Let $m=2k$. Then $P^{(1)}(z)=bz^k + Q^{(1)}(z)$, $b\ne 0$ and
$\deg Q^{(1)}\leq k-1$. Define $q(z):=p(z)-be_{2k}=
\begin{pmatrix} Q^{(1)}(z)\\ Q^{(2)}(z)
\end{pmatrix}$.
Here $\deg Q^{(2)}\leq k-1$, so $h(q)\leq 2k-1=m-1$ and we can again
apply the inductive hypothesis. \hfill $\square$

\medskip

The latter proposition can be extended in a natural way.

\begin{proposition} \label{propos2.3}
{\it Let $\{g_n(z)\}_{n\ge0}$ be an arbitrary sequence of
vector-polynomials, such that
$$ h(g_n) = n, \qquad  n\in\dZ_+.  $$
Then $\{g_n\}_{n\ge 0}$ is a basis in $\dP_\infty$. }
\end{proposition}

{\it Proof.} By Proposition \ref{propos2.2}
$$ g_m(z) = \sum_{k=0}^m c_{m,k} e_k(z), \quad c_{m,m}\ne 0,
$$
or in a vector-matrix form
$$ \begin{pmatrix} g_0 \\ \vdots \\ g_m
\end{pmatrix} =
\begin{pmatrix} c_{00} &  &  & \\
c_{10} & c_{11} & & \\
\vdots & & \ddots & \\
c_{m0} & \ldots & \ldots & c_{mm}
\end{pmatrix}
\begin{pmatrix} e_0 \\ \vdots \\ e_m
\end{pmatrix},
$$
where the matrix of $(c_{m,k})$ is triangular. So,
$$ \begin{pmatrix} e_0 \\ \vdots \\ e_m
\end{pmatrix} =
\begin{pmatrix} \wt c_{00} &  & & \\
\wt c_{10} & \wt c_{11} & &  \\
\vdots & & \ddots & \\
\wt c_{m0} & \ldots & \ldots & \wt c_{mm}
\end{pmatrix}
\begin{pmatrix} g_0 \\ \vdots \\ g_m
\end{pmatrix}, \qquad \wt c_{m,m} \ne 0.
$$
Since $\{e_n\}_{n\ge 0}$ is a basis, then so is $\{g_n\}_{n\ge 0}$.
\hfill $\square$
\medskip

\subsection{Transforms in $\dP_\infty$}

A $2\times 2$ matrix $ A = \begin{pmatrix} a & b \\ c & d
\end{pmatrix}$ defines a transform of the vector-polynomials:
$$ A p(z) = A
\begin{pmatrix} P^{(1)}(z)\\ P^{(2)}(z)
\end{pmatrix} =
\begin{pmatrix} a P^{(1)}(z) + b P^{(2)}(z)\\
c P^{(1)}(z) + d P^{(2)}(z)
\end{pmatrix}
$$

\begin{proposition} \label{propos2.4}
{\it
\begin{enumerate}
    \item  For arbitrary matrix $A$ and $p\in \dP_\infty$
    $$ h(Ap) \leq h(p) + 1. $$
    \item  If $A$ is upper-triangular, then for arbitrary $p\in\dP_\infty$
    $$ h(Ap)\leq h(p). $$
    \item  If $A$ is lower-triangular, then
    $$h(p)\le 2k+1 \Longrightarrow h(Ap)\le 2k+1. $$
\end{enumerate}
}
\end{proposition}

{\it Proof.} \ (1) By (\ref{2.6})
$$ h(Ap) = \max \(2\deg(aP^{(1)}+bP^{(2)}), \ 2\deg (c P^{(1)} + d
P^{(2)})+1\),$$

$$ \begin{aligned} 2\deg (aP^{(1)}+bP^{(2)}) &\leq 2\max (\deg P^{(1)},\ \deg P^{(2)})
=\max (2\deg P^{(1)}, \ 2\deg P^{(2)}) \leq h(p), \\
2\deg (cP^{(1)} + d P^{(2)} ) &\leq \max (2\deg P^{(1)}, \ 2\deg
P^{(2)}) \leq h(p), \end{aligned} $$
$$ 2\deg (cP^{(1)}+dP^{(2)})+1\leq h(p)+1.
$$

(2) If $c=0$, then $2\deg (dP^{(2)})+1\leq 2\deg P^{(2)} +1 \leq
h(P)$.

(3) Let now $h(p)\le 2k+1$ and $ A = \begin{pmatrix} a & 0 \\ c & d
\end{pmatrix}$. Then
$$ A p(z) = \begin{pmatrix} a P^{(1)}(z)\\ c P^{(1)} (z) + d P^{(2)}(z)
\end{pmatrix},
$$
and we have by assumption $2\deg P^{(1)} \leq 2k+1$,\  $2\deg
P^{(2)}\leq 2k$. So,
$$\deg(cP^{(1)} + dP^{(2)})\leq k \Longrightarrow h(Ap)\leq 2k+1, $$
as claimed.   \hfill $\square$

\bigskip

Later on we will use the following property of the height, which is
a simple consequence of (\ref{2.6}): for
$p=\begin{pmatrix} P^{(1)}\\
P^{(2)} \end{pmatrix} \in \dP_\infty$ we have

\begin{equation} \label{2.61}
\left\{ \begin{array}{l} h(p)\leq 2k+1 \Longrightarrow h
\begin{pmatrix} (z-a) P^{(1)}(z)\\ P^{(2)}(z)
\end{pmatrix}\leq 2k+2,\\
h(p)\leq 2k \Longrightarrow h
\begin{pmatrix} P^{(1)}(z)\\ (z-a) P^{(2)}(z)
\end{pmatrix}\leq 2k+1.\\
\end{array}
\right. \end{equation}

\subsection{The generators of interpolation problem}

It is clear that solutions of (\ref{2.3}) form a module over the
ring of polynomials in $\dP_\infty$, i.e., if $r$ and $q$ are
solutions of \eqref{2.3}, then so is $Sr+Tq$ for arbitrary
polynomials $S$ and $T$. The goal of this subsection is to show that
this module has exactly two generators, and to study their
properties.

Recall that $\dP_\infty(I_n)$ is the set of all solutions of
(\ref{2.3}). Set
$$ h(I_n) := \min \bigl\{ h(q): \ q\in\dP_\infty(I_n), \quad q\ne 0 \bigr\},
$$
which we call the height of the interpolation problem.

\begin{definition}\label{def2}
We say that $r=\in \dP_\infty(I_n)$ is a minimal generator of
$(\ref{2.3})$, if
$$ h(r) = h(I_n).
$$
\end{definition}

\begin{proposition} \label{propos2.8}
{\it The minimal generator of $(I_n)$ is unique up to a constant
factor.}
\end{proposition}

{\it Proof.} \ Let $r_1$ and $r_2$ be two minimal generators. By
Proposition \ref{propos2.1}, $\exists a \in \dC$: $h(ar_1 + r_2)
\leq h(I_n) -1$. But $ar_1+br_2$ is a solution of (\ref{2.3}). Since
$r_1$ and $r_2$ are minimal non-trivial solutions of (\ref{2.3}), we
conclude $a r_1 + b r_2=0$. \hfill $\square$

\medskip

A trivial (nonzero) solution of \eqref{2.3}
$P^{(1)}=P^{(2)}=\prod_j(z-z_j)$ provides the bound $h(I_n)\le
2n+1$. It turns out that this bound can be improved immensely.

\begin{theorem}\label{theor2.6}
$ h(I_n) \leq n$.
\end{theorem}

{\it Proof.} \ The following non-negative matrices of rank 1 play a
key role in our consideration:
$$ \s_k =
\begin{pmatrix} |\a_{k,1}|^2 & \bar \a_{k,1} \a_{k,2}\\
\bar \a_{k,2} \a_{k,1} & |\a_{k,2}|^2
\end{pmatrix} =
\begin{pmatrix} \bar\a_{k,1} \\
\bar\a_{k,2}
\end{pmatrix}
(\a_{k,1}, \,\a_{k,2}), \qquad k=1,\ldots,n.
$$
It is clear that
$$ |\a_{k,1} P^{(1)}(z_k) + \a_{k,2} P^{(2)}(z_k)|^2=
p^*(z_k) \s_k\, {p(z_k)},
$$
$$ p^*(z) :=
(\ov{P^{(1)}(z)},\, \ov{P^{(2)}(z)}), \qquad p(z):=
\begin{pmatrix} P^{(1)}(z) \\ P^{(2)}(z)
\end{pmatrix},
$$
so the problem (\ref{2.3}) is equivalent to
\begin{equation} \label{2.7}
p^*(z_k) \s_k\, p(z_k) = 0, \qquad k=1,\ldots,n.
\end{equation}

We proceed by induction on $n$.

1. For $n=1$ we have a non-trivial solution
$$ p:= \begin{pmatrix} -\a_{1,2} \\ \a_{1,1}
\end{pmatrix}, \qquad h(p)\leq 1.
$$

2. Suppose that we have already proved the result for $n$, and we
want to prove it for $n+1$. The forthcoming construction depends on
whether $n$ is odd or even.

Let $n=2k+1$, and consider the problem $(I_{n+1})$ \eqref {2.7} with
$n+1$ data. If $\a_{j,1}=0$ for all $j=1,\ldots,n+1$, then the
vector-polynomial $p=
\begin{pmatrix} 1 \\ 0 \end{pmatrix}$, \ $h(p)=0$, is a solution of
$(I_n)$, and we are done. So, suppose without loss of generality
that $\a_{n+1,1}\ne 0$ (otherwise enumerate the points
$z_1,\ldots,z_{n+1}$). The upper-triangular matrix
$$ \O_0 := \begin{pmatrix} (\a_{n+1,1})^{-1} &
- \a_{n+1,2} (\a_{n+1,1})^{-1} \\
0 & 1 \end{pmatrix},
$$
satisfies
\begin{equation}\label{2.71}
\O_0^*\, \s_{n+1}\, \O_0 =
\begin{pmatrix} 1 & 0 \\
0 & 0 \end{pmatrix}, \end{equation}
 and
\begin{equation} \label{2.8}
 \O_0^*\, \s_j\, \O_0 =
\begin{pmatrix} |\b_{j,1}|^2 & \bar \b_{j,1} \b_{j,2}\\
\bar \b_{j,2} \b_{j,1} & |\b_{j,2}|^2
\end{pmatrix};
\qquad j=1,\ldots,n,
\end{equation}
with some numbers $\b_{j,1}, \b_{j,2}, \ j=1,\ldots,n$. Put
$$ \g_{j,1}:= (z_{n+1} - z_j) \b_{j,1}; \qquad \g_{j,2}
:= \b_{j,2}, \qquad j=1,\ldots,n,
$$
and consider an auxiliary interpolation problem $(\wt I_n)$:
$$ \g_{j,1} P^{(1)} (z_j) + \g_{j,2} P^{(2)} (z_j) = 0,
\qquad j=1,\ldots,n.
$$
By the induction hypothesis, there exists a solution
$q\in\dP_\infty$, \ $h(q) \leq n=2k+1$, of this problem:
$$ q^*(z_j)\wt\s_j\, q(z_j) =0, \qquad j=1,\ldots,n, $$
where
$$\wt\s_j =\begin{pmatrix} |\g_{j,1}|^2 & \bar\g_{j,1} \g_{j,2} \\
\bar\g_{j,2} \g_{j,1} & |\g_{j,2}|^2
\end{pmatrix}
=\begin{pmatrix} \ov{z_{n+1}-z_j} & 0 \\
0 & 1 \end{pmatrix} \O_0^* \s_j \O_0
\begin{pmatrix} z_{n+1}-z_j & 0 \\
0 & 1 \end{pmatrix}.
$$

Define a vector-polynomial
\begin{equation}\label{2.9} r(z) =
\begin{pmatrix} R^{(1)}(z) \\ R^{(2)}(z)
\end{pmatrix} := \O_0
\begin{pmatrix} z_{n+1}-z & 0 \\
0 & 1 \end{pmatrix} q(z).
\end{equation}
For $r$ we have
$$ r^*(z_j)\, \s_j \, r(z_j) = q^*(z_j)\, \wt\s_j \, q(z_j) = 0, \quad
j=1,\ldots,n,
$$
and for $j=n+1$
$$ r^*(z_{n+1})\, \s_{n+1} \, r(z_{n+1}) = q^*(z_{n+1})
\begin{pmatrix} 0 & 0 \\ 0 & 1
\end{pmatrix}
\begin{pmatrix} 1 & 0 \\ 0 & 0
\end{pmatrix}
\begin{pmatrix} 0 & 0 \\ 0 & 1
\end{pmatrix}
q(z_{n+1}) = 0,
$$
by \eqref{2.71}. So, $r$ is a solution of $(I_{n+1})$.

Since $h(q)\leq 2k+1$, then by the upper inequality in \eqref{2.61}
we have
$$
h \begin{pmatrix} (z_{n+1}-z) Q^{(1)}(z) \\
Q^{(2)}(z)
\end{pmatrix} \leq 2k+2 = n+1,
$$
so, by Proposition \ref{propos2.4}, $h(r)\leq n+1 $, as needed.

\smallskip
Let $n=2k$, and $(I_{n+1})$ \eqref {2.7} be the
interpolation problem with $n+1$ data. If $\a_{j,2}=0$ for all
$j=1,\ldots,n+1$, then the vector-polynomial $p=
\begin{pmatrix} 0 \\ 1 \end{pmatrix}$, \ $h(p)=1$, is a solution of
$(I_n)$, and we are done. So, suppose as above, that $\a_{n+1,2}\ne
0$. The lower-triangular matrix
$$ \O_1 := \begin{pmatrix} 1 & 0 \\
- \a_{n+1,1} (\a_{n+1,2})^{-1} & (\a_{n+1,2})^{-1} \end{pmatrix},
$$
satisfies
\begin{equation}\label{2.10}
\O_1^*\, \s_{n+1}\, \O_1 =
\begin{pmatrix} 0 & 0 \\ 0 & 1 \end{pmatrix}, \end{equation}
 and
\begin{equation} \label{2.101}
 \O_1^*\, \s_j\, \O_1 =
\begin{pmatrix} |\d_{j,1}|^2 & \bar \d_{j,1} \d_{j,2}\\
\bar \d_{j,2} \d_{j,1} & |\d_{j,2}|^2
\end{pmatrix};
\qquad j=1,\ldots,n.
\end{equation}
Put
$$ \l_{j,2}:= (z_{n+1} - z_j) \d_{j,2}; \qquad \l_{j,1}
:= \d_{j,1}, \qquad j=1,\ldots,n,
$$
and consider an auxiliary interpolation problem $(I_n)$:
$$ \l_{j,1} P^{(1)} (z_j) + \l_{j,2} P^{(2)} (z_j) = 0,
\qquad j=1,\ldots,n,
$$
By the induction hypothesis, there exists a solution
$q\in\dP_\infty$, \ $h(q) \leq n=2k$, of this problem:
$$ q^*(z_j)\wt\s_j\, q(z_j) =0, \qquad j=1,\ldots,n, $$
where
$$\wt\s_j =\begin{pmatrix} |\l_{j,1}|^2 & \bar\l_{j,1} \l_{j,2} \\
\bar\l_{j,2} \l_{j,1} & |\l_{j,2}|^2
\end{pmatrix}
=\begin{pmatrix} 1 & 0 \\
0 & \ov{z_{n+1}-z_j} \end{pmatrix} \O_1^* \s_j \O_1
\begin{pmatrix} 1 & 0 \\
0 & z_{n+1}-z_j \end{pmatrix}.
$$

Define a vector-polynomial
\begin{equation}\label{2.102} r(z) =
\begin{pmatrix} R^{(1)}(z) \\ R^{(2)}(z)
\end{pmatrix} := \O_1
\begin{pmatrix} 1 & 0 \\
0 & z_{n+1}-z \end{pmatrix} q(z).
\end{equation}
For $r$ we have
$$ r^*(z_j)\, \s_j \, r(z_j) = q^*(z_j)\, \wt\s_j \, q(z_j) = 0, \quad
j=1,\ldots,n,
$$
and for $j=n+1$
$$ r^*(z_{n+1})\, \s_{n+1} \, r(z_{n+1}) = q^*(z_{n+1})
\begin{pmatrix} 1 & 0 \\ 0 & 0
\end{pmatrix}
\begin{pmatrix} 0 & 0 \\ 0 & 1
\end{pmatrix}
\begin{pmatrix} 1 & 0 \\ 0 & 0
\end{pmatrix}
q(z_{n+1}) = 0,
$$
by \eqref{2.10}. So, $r$ is a solution of $(I_{n+1})$.

Since $h(q)\leq 2k$, then by the lower inequality in \eqref{2.61} we
have
$$
h \begin{pmatrix}  Q^{(1)}(z) \\
(z_{n+1}-z)Q^{(2)}(z)
\end{pmatrix} \leq 2k+1 = n+1,
$$
so, by (3), Proposition \ref{propos2.4} ($\O_1$ is
lower-triangular), $h(r)\leq 2k+1=n+1 $. The proof is complete.
\hfill $\square$

\bigskip

If $r$ is a minimal generator of $(I_n)$, then $Sr\in
\dP_\infty(I_n)$ for any polynomial $S$. So, the question arises
naturally whether $\dP_\infty(I_n)=\{Sr\}$, $S$ a polynomial. The
answer is negative: it turns out that the module of solutions of
$(I_n)$ has exactly one more generator. Denote
$\dP'_\infty(I_n)=\dP_\infty(I_n)\backslash\{Sr\}$. It is shown in
Theorem \ref{theor2.7} below that this set is nonempty, so the
following definition makes sense.

\begin{definition}\label{def3} We say that $q \in
\dP_\infty(I_n)$ is a second generator of $(I_n)$, if
$$ h(q) = min \{ h(p), \ p\in \dP'_\infty(I_n)\}.
$$
\end{definition}

\begin{theorem}\label{theor2.7}
The set $\dP'_\infty(I_n)$ is nonempty. Furthermore, the height of
any second generator is $h(q)=2n+1-h(I_n)$.
\end{theorem}

{\it Proof.} We show that there is a solution of $(I_n)$ of the
height $\le 2n+1-h(I_n)$, which is not of the form $Sr$.

Let $h(I_n)=k\le n$. Pick arbitrary different numbers
$z_{n+1},z_{n+2},\ldots,z_{2n+1-k}\in \dC$ distinct from
$z_1,\ldots,z_n$. Take numbers $\a_{j,1}$, $\a_{j,2}$,
$j=n+1,n+2,\ldots,2n+1-k$ in such a way that
\begin{equation}\label{2.103} \a_{j,1}
R^{(1)} (z_j) + \a_{j,2} R^{(2)} (z_j) \ne 0, \quad
j=n+1,n+2,\ldots,2n+1-k,
\end{equation}
where $ r = \left(\begin{array}{c} R^{(1)} \\ R^{(2)}\\
\end{array}\right) $ is the minimal generator of $(I_n)$.
Consider the interpolation problem
\begin{equation}\label{2.11} \a_{j,1}
P^{(1)} (z_j) + \a_{j,2} P^{(2)} (z_j) = 0, \quad
j=1,2,\ldots,2n+1-k.
\end{equation}
By Theorem \ref{theor2.6}, there exists a nonzero vector-polynomial
$ p =
 \left(\begin{array}{c} P^{(1)} \\ P^{(2)}\\
\end{array}\right) \ne 0$, which solves this problem, and $h(p)
\leq 2n + 1 -k$.

Suppose that $p = Sr$ for a polynomial $S$. Then, by (\ref{2.11}),
$$ S(z_j) \bigl[ \a_{j,1} R^{(1)} (z_j) + \a_{j,2} R^{(2)} (z_j) \bigr] = 0,
\quad j=n+1,n+2,\ldots,2n+1-k,
$$
so by (\ref{2.103}),
$$ S(z_{n+1}) = S(z_{n+2}) = \ldots = S(z_{2n+1-k})=0.
$$
Since $S(z)\not\equiv 0$, we see that $\deg S(z) \ge n+1 -k$, and by
\eqref{2.61}
$$ h(p)=h(Sr)=h(r)+2\deg S\ge h(r) + 2(n+1-k) = 2n +2 -k,
$$
which leads to contradiction with $h(p) \leq 2n +1 -k$. So
$p\in\dP'_\infty(I_n)$ and $h(p)\le 2n+1-k$, as claimed.

\medskip

Let now $q$ be any second generator, so $h(q)\le 2n+1-k$. We prove
next $h(q) \ge 2n+1-k$.

We have
$$
\left\{ \begin{array}{l} \a_{j,1} R^{(1)} (z_j) + \a_{j,2} R^{(2)} (z_j) =0,\\
\a_{j,1} Q^{(1)} (z_j) + \a_{j,2} Q^{(2)} (z_j) =0 ,\\
\end{array}
\right. \qquad j=1,\ldots,n,
$$
so
\begin{equation} \label{2.12}
\det \begin{pmatrix} R^{(1)}(z_j) & R^{(2)}(z_j) \\ Q^{(1)}(z_j) &
Q^{(2)}(z_j) \end{pmatrix}=
  R^{(1)} (z_j) Q^{(2)} (z_j) - R^{(2)} (z_j) Q^{(1)} (z_j)=0,
\end{equation}
$j=1,\ldots,n$. Suppose that $h(q) \le 2n-k$, which implies by
\eqref{2.6}
$$\deg Q^{(1)} \leq n - \frac{k}{2}, \qquad
\deg Q^{(2)} \leq n - \frac{k}{2} - 1.
$$
If $h(r)=k$ is even, then
$$ \deg R^{(1)} = \frac{k}{2}, \quad \deg R^{(2)} \leq \frac{k}{2}
-2,
$$
and so
$$ \deg \bigl( R^{(1)} Q^{(2)} - R^{(2)} Q^{(1)} \bigr) \leq n-1.
$$
If $k$ is odd, then
$$ \deg R^{(2)} = \frac{k-1}{2}, \quad \deg R^{(1)} \leq \frac{k-1}{2},
$$
and again
$$ \deg \bigl( R^{(1)} Q^{(2)} - R^{(2)} Q^{(1)} \bigr) \leq n-1.
$$
It follows now from \eqref{2.12} that
 \begin{equation} \label{2.121}
 R^{(1)}(z)Q^{(2)}(z)-R^{(2)}(z)Q^{(1)}(z) \equiv 0.
\end{equation}

Let $T$ be the greatest common divisor of $R^{(1)}$ and $R^{(2)}$,
so (\ref{2.121}) turns into $X^{(1)}TQ^{(2)}=X^{(2)}TQ^{(1)}$, with
relatively prime $X^{(1)}:=R^{(1)}/T$ and $X^{(2)}:=R^{(2)}/T$.
Hence $S^{(1)}:=Q^{(1)}/X^{(1)}$ is a polynomial, and so is
$S^{(2)}:=Q^{(2)}/X^{(2)}$. Now $X^{(1)}TQ^{(2)}=X^{(2)}TQ^{(1)}$
implies $S^{(1)}=S^{(2)}=S$ and
\begin{equation} \label{2.122}
Q^{(i)}=\frac{R^{(i)}}{T} \cdot S\,, \qquad i=1,2.
\end{equation}
Note that all roots of $T$ are among the nodes of interpolation.
Indeed, if $T(w)=0$ and $w\not\in\{z_1,\ldots,z_n\}$, then
$$ \wh r(z)= \begin{pmatrix} \wh R^{(1)} \\ \wh R^{(2)}
\end{pmatrix}= \frac1{z-w}\,\begin{pmatrix}  R^{(1)} \\ R^{(2)}
\end{pmatrix} \in\dP_\infty(I_n) $$
and $h(\wh r)<h(r)$, which is impossible since $r$ is a minimal
generator of $(I_n)$. Hence, $w=z_l$.

It remains only to show that $S(z_l)=0$. Assume that $S(z_l) \ne 0$.
Then
$$ \a_{l,1}Q^{(1)}(z_l)+\a_{l,2}Q^{(2)}(z_l)=0 $$
and \eqref{2.122} imply
$$ \lim_{z\to z_l} \bigl\{ \a_{l,1}\,\frac{R^{(1)}(z)}{T(z)}+
\a_{l,2}\,\frac{R^{(2)}(z)}{T(z)} \bigr\}=0. $$ So
$$ \lim_{z\to z_l} \bigl\{ \a_{l,1}\,\frac{R^{(1)}(z)}{(z-z_l)^{n_l}}+
\a_{l,2}\,\frac{R^{(2)}(z)}{(z-z_l)^{n_l}} \bigr\} =0, $$ where
$T=(z-z_l)^{n_l}\wt T$, $\wt T(z_l)\ne 0$. Since
$R^{(i)}=(z-z_l)^{n_l}\wt R^{(i)}$, $i=1,2$, we see that
$$ \wt r(z)=\frac1{(z-z_l)^{n_l}}\, r(z)\in\dP_\infty(I_n) $$
and clearly $h(\wt r)<h(r)$, which again leads to contradiction with
$r$ being the minimal generator.

Finally, since all the roots of $T$ are among $\{z_1,\ldots,\z_k\}$
and $S(z_1)=\ldots=S(z_k)=0$, $P=S/T$ is a polynomial, and by
\eqref{2.122} $Q^{(i)}=PR^{(i)}$, which contradicts to
$q\in\dP'_\infty(I_n)$. The proof is complete. \hfill $\square$

{\it  Remark}. In fact we have proven that each solution
$q\in\dP'_\infty(I_n)$ with $h(q)\le 2n+1-h(I_n)$ is a second
generator. It is also not hard to see that each solution
$q\in\dP(I_n)$ with $h(q)= 2n+1-h(I_n)$ is a second generator.

\begin{theorem}\label{theor2.8}
Each solution of the problem $(I_n)$ has the form
\begin{equation} \label{2.14}
p(z) = S(z) r(z) + T(z) q(z),
\end{equation}
where $r$ and $q$ are the minimal and second generators of $(I_n)$,
respectively, and $S$ and $T$ are polynomials. Conversely, each
vector-polynomial of the form $(\ref{2.14})$ with arbitrary
polynomials $S$ and $T$ belongs to $\dP_\infty(I_n)$.
\end{theorem}

{\it Proof. } We only prove the first statement. Let $h(I_n)=k$.
Consider a system of vector-polynomials $\{f_j\}_{j\ge 0}$ defined
as follows:
$$ \begin{aligned} f_k &=r, \qquad f_{2n+1-k}=q, \\
f_j &=e_j, \quad j=0,1,\ldots,k-1; \qquad f_{k+2j-1}=e_{k+2j-1},
\quad j=1,2,\ldots,n-k, \\
f_{k+2j} &=z^j r,  \quad f_{2n+1-k+2j}=z^j q; \qquad j\in\dN,
\end{aligned} $$
$e_i$ are in \eqref{2.62}. It is easy to check that $h(f_j)=j$ for
all $j$, so by Proposition \ref{propos2.3} this system is a basis in
$\dP_\infty$, and in particular each $p\in\dP_\infty(I_n)$ admits a
unique representation in the form
$$
p(z) = \sum_{i=0}^{k-1} a_i e_i + \sum_{i=0}^{n-1-k} b_i e_{k+1+2i}
+ S_1(z) r(z) + T_1(z) q(z).
$$
Since $p, r, q\in\dP_\infty(I_n)$ then $\displaystyle{\sum_0^{k-1}
a_i e_i + \sum_{i=0}^{n-1-k} b_i e_{k+1+2i}}$ is also solution of
$(I_n)$. But its height is less then $2n+1-k$, so, according to the
definition of a second generator,
$$ \sum_0^{k-1} a_i e_i + \sum_{i=0}^{n-1-k} b_i e_{k+1+2i} =
\wt S(z) r(z).
$$
Hence $p$ is of the form \eqref{2.14} with $S=S_1+\wt S$, $T=T_1$.
\hfill $\square$

\medskip

Thus, the following theorem holds for interpolation problem
(\ref{2.1}):

\begin{theorem}\label{theor2.9}
Each solution of problem $(\ref{2.1})$ has the form
\begin{equation} \label{2.16}
\frac{ S(z) R^{(1)}(z) + T(z) Q^{(1)}(z)}{S(z) R^{(2)}(z) + T(z)
Q^{(2)}(z)}\,,
\end{equation}
where $r=\begin{pmatrix} R^{(1)} \\ R^{(2)}
\end{pmatrix}$ and $q=\begin{pmatrix} Q^{(1)} \\ Q^{(2)}
\end{pmatrix}$ are minimal and second generators of $(I_n)$,
and $S$ and $T$ are polynomials. Conversely, if $r$ and $q$ are the
minimal and  second generators of $(I_n)$, and $S$ and $T$ are such
polynomials that the numerator and the denominator in $(\ref{2.16})$
have no common roots, then $(\ref{2.16})$ is a solution of
$(\ref{2.1})$.
\end{theorem}

\bigskip

{\it Remark.} Roughly speaking, the task of giving {\it a}
description for the set of the solutions of rational interpolation
problem (\ref{2.1}), is obvious. Let
${\ds\frac{R^{(1)}}{R^{(2)}}}$ be {\it any} rational function
solving problem (\ref{2.1}). Then all the functions of the type
${\ds\frac{R^{(1)}}{R^{(2)}} + \frac{Q^{(1)}}{Q^{(2)}}}$, where
${\ds\frac{Q^{(1)}}{Q^{(2)}}}$ is an arbitrary rational function
vanishing in the nodes of interpolation, will be all the solutions
of the rational interpolation problem. However, such a description
does not permit us to predict the degrees of the numerator and
denominator as well as other properties needed in applications.
So, we cannot obtain in this way a rational function which solves
the interpolation problem and has the prescribed properties. For
example, in the interpolation problem appearing in the next
section we will need a rational function with monic numerator and
denominator of concrete degrees, such that the free term of the
numerator equals 1. Thus, more elaborated results are required.
Certainly, we do not think that the two descriptions for the
solutions of the interpolation problem, mentioned above, are the
only possible.

\setcounter{section}{2} \setcounter{equation}{0}
\section {Reduction of MISP to Rational Interpolation}

For the definitions, notations and basic properties of finite CMV
matrices see, for example, \cite{sim2, GK2}. We will add some more
to the list.

Let $\cC=\cC(\alpha_0,\ldots,\alpha_{n-2};\beta)$ be a finite CMV
matrix with Verblunsky's parameters $(\a_0,\ldots,\a_{n-2};\b)$ and
the system of the Szeg\H{o} polynomials
$\{\F_0,\ldots,\F_{n-1};\wt\F_n\}$. They satisfy the famous {\it
\SZ\ recurrence relations}
\begin{equation}\begin{aligned}
\label{1.0} \Phi_k(z) &=z\Phi_{k-1}(z) - \bar\alpha_{k-1}
\Phi_{k-1}^*(z), \qquad
k=1,2,\ldots,n-1, \quad \Phi_0\equiv 1, \\
\wt\Phi_n(z) &=z\Phi_{n-1}(z) - \bar\beta \Phi_{n-1}^*(z).
\end{aligned} \end{equation}
As is known,
$$\wt\F_n(z)=\prod_{j=1}^n(z-\z_j), \qquad \S(\cC)=\{\z_j\}_1^n $$ a
spectrum of $\cC$, $\z_j\ne\z_i$, $j\ne i$. Put
\begin{equation} \label{1.1}
\k_m=\prod_{j=0}^{m-1} (1-|\a_j|^2)^{-1/2}, \quad m=0,1,\ldots,n-1;
\qquad \k_0=1,
\end{equation}
and for appropriate values of the indices define
\begin{equation} \label{1.2}
x_{2k}(z):=z^{-k}\k_{2k}\F_{2k}(z), \quad
x_{2k+1}(z):=z^{-k-1}\k_{2k+1}\F_{2k+1}^*(z),
\end{equation}
where by definition $\varphi_k^*(z) = z^k \overline{\varphi_k
(1/\overline{z})}$.

For each eigenvalue $\z_j$ the following equality
\begin{equation} \label{1.3}
\cC X_j=\z_jX_j, \qquad X_j=[x_0(\z_j),\ldots,x_{n-1}(\z_j)]^t ,
\end{equation}
gives (along with \eqref{1.2}) an explicit expression for the
eigenvectors of $\cC$ in terms of the Szeg\H{o} polynomials and
Verblunsky parameters. \eqref{1.3} is proved in \cite[Lemma
4.3.14]{simA}, for infinite CMV matrices. For finite matrices the
argument is similar. As a matter of fact, the following more precise
result holds.

\begin{proposition} \label{propos1.1}
For $z\in\dC\backslash\{0\}$ and $X(z)=[x_0(z),\ldots,x_{n-1}(z)]^t$
the equality holds
$$ (z-\cC)X=z^{-[n/2]}\k_{n-1}\wt\F_n(z)v_n, \qquad v_n\in\dC^n,
\quad \|v_n\|=1. $$
\end{proposition}

Due to the sieving procedure we will assume without loss of
generality that $n$ is an even number: $n=2l$. Throughout the rest
of the paper we assume that the last Verblunsky coefficient is
known, and for simplicity put $\b=1$. The $\cL\cM$ factorization now
takes the form $\cC(\a_0,\ldots,\a_{2l-2};1)=\cL\cM$ with
\begin{equation} \label{1.4}
\cL=\begin{pmatrix} \te(\a_0) & \ & \ \\
                     \  & \ddots & \ \\
                     \ & \ & \te(\a_{2l-2})
\end{pmatrix}, \qquad
\cM=\begin{pmatrix} 1 & \ & \ & \ & \ \\
                    \ & \te(\a_1) & \ & \ & \ \\
                    \ & \ & \ddots & \ & \ \\
                    \ & \ & \ & \te(\a_{2l-3}) & \ \\
                    \ & \ & \ & \ & \ & 1
\end{pmatrix},
\end{equation}
$$
\Theta(\a_j)=\begin{pmatrix} \bar{\alpha_j} & \rho_j \\
\rho_j & -\alpha_j \\
\end{pmatrix}, \qquad |\alpha_j|<1, \quad \rho_j=\sqrt{1-|\a_j|^2}>0.
$$

Put
\begin{equation} \label{1.5}
U =\begin{pmatrix}  \dO & \ldots & \dO & J  \\
                    \dO & \ldots & J & \dO  \\
                    \ldots & \ldots & \ldots & \ldots \\
                    J & \dO & \ldots & \dO
\end{pmatrix},  \qquad
J = \begin{pmatrix} 0 & 1 \\
                    1 & 0
\end{pmatrix},
\end{equation}
the orthogonal $2l\times 2l$ matrix, and consider the {\it
reflection} of $\cC$
\begin{equation} \label{1.6}
\cC_r := U \cC U = U \cL U \cdot U \cM U = \cL_r \cdot \cM_r.
\end{equation}
It is clear from (\ref{1.4})--\eqref{1.6}, that
\begin{equation} \label{1.7}
\begin{aligned}
\cL_r &= \begin{pmatrix} \te(-\bar\a_{2l-2}) & \ & \ \\
                     \  & \ddots & \ \\
                     \ & \ & \te(-\bar\a_0)
\end{pmatrix}; \\
\cM_r &= \begin{pmatrix} 1 & \ & \ & \ & \ \\
                    \ & \te(-\bar\a_{2l-3}) & \ & \ & \ \\
                    \ & \ & \ddots & \ & \ \\
                    \ & \ & \ & \ & \te(-\bar\a_1)\\
                    \ & \ & \ & \ & \ & 1
\end{pmatrix},
\end{aligned}
\end{equation} so
\begin{equation} \label{1.71} \cC_r=\cC (\l_0,\ldots,\l_{2l-2};1), \qquad \l_k := -
\bar\a_{2l-2-k}, \ k=0,1,\ldots,2l-2, \end{equation} is also a CMV
matrix corresponding to the ``reversed'' Verblunsky parameters. We
denote the Szeg\H{o} polynomials for $\cC_r$ by
$\{\L_{0},\ldots,\L_{n-1};\wt\L_{n}\}$.

Obviously, $\S(\cC_r)=\S(C)=\{\z_j\}_1^n$ and
$$ \cC_r Y_{j} = \z_j Y_{j}, \qquad Y_{j} =
[y_{0}(\z_j),\ldots, y_{n-1}(\z_j)]^t,
$$
where $y_{k}$ are in (\ref{1.2}) for the matrix $\cC_r$. On the
other hand, by (\ref{1.3}) and (\ref{1.6})
$$ \cC_r\hat X_j = \z_j \hat X_j, \quad \hat X_j =
[x_{n-1}(\z_j),\ldots,x_0(\z_j)]^t.
$$
Since the spectrum is simple, the vectors $Y_{j}$ and $\hat X_j$ are
proportional:
$$ Y_{j} = c_j \hat X_j, \quad y_{k-1}(\z_j)=c_j x_{n-k}(\z_j);
\quad k,j = 1,2,\ldots,n,
$$
or
\begin{equation} \label{1.8}
\frac{y_{k-1}(\z_j)}{y_{k}(\z_j)} =
\frac{x_{n-k}(\z_j)}{x_{n-k-1}(\z_j)}, \qquad k=1,2,\ldots,n-1.
\end{equation}

\bigskip

Under the {\it mixed inverse spectral problem} (MISP) we mean the
reconstrucion of a CMV matrix $\cC=\cC(\a_0,\ldots,\a_{n-2};1)$, or
equivalently, of a set of Verblunsky parameters
$\a_0,\ldots,\a_{n-2} \in \dD$, when a part $\{\z_j\}_{j=1}^m$ of
its spectrum and a part of the system $(\a_0,\ldots,\a_{n-2})$ are
known.


Here is the simplest problem of this type. Assume that we know
$(\a_0,\ldots,\a_{n-3})$ as well as two eigenvalues $\z_1\ne\z_2$,
and $\a_{n-2}$ is to be found so that $\z_{1,2}\in\S(\cC)$. Once
$\F_{n-2}$ is known, we apply the Szeg\H{o} recurrences to obtain
$$ \wt\F_n(z) = z (z+\a_{n-2}) \F_{n-2}(z) - (z \bar\a_{n-2} + 1)
\F_{n-2}^*,
$$
so
\begin{equation}\label{1.9}
b(\z_j) = \tau_j, \quad b(\l) =
\frac{\l+\a_{n-2}}{1+\l\bar\a_{n-2}}\,, \qquad j=1,2,
\end{equation}
$$ \tau_j=\frac{\Phi_{n-2}^*(\z_j)}{\z_j \F_{n-2}(\z_j)}\,,\qquad
j=1,2. $$
The question is whether $\a_{n-2}$ is uniquely
determined from the interpolation problem (\ref{1.9}). An
elementary analysis of (\ref{1.9}) shows that it has a unique
solution as long as $\tau_1\z_1\ne\tau_2\z_2$, that is,
\begin{equation}\label{1.10}
\frac{\F_{n-2}^*(\z_1)}{\F_{n-2}(\z_1)} \ne
\frac{\F_{n-2}^*(\z_2)}{\F_{n-2}(\z_2)}\,,
\end{equation}
it has infinitely many solutions if $\t_2=-\z_1$ and $\t_1=-\z_2$,
or
$$ \frac{\F_{n-2}^*(\z_1)}{\F_{n-2}(\z_1)} =
\frac{\F_{n-2}^*(\z_2)}{\F_{n-2}(\z_2)} = -\z_1\z_2\,,
$$
and it has no solutions at all, if $\tau_1\z_1 = \tau_2\z_2$, but
$\tau_1 \ne -\z_2$, $\tau_2 \ne -\z_1$. It is not hard to check that
each situation may occur for interpolation problem \eqref{1.9}.
However, if the existence of CMV matrix $\cC$ with
$\z_{1,2}\in\S(\cC)$ is supposed, the existence of the solution of
problem \eqref{1.9} is guaranteed and the problem of finding
$\a_{n-2}$ may have either unique or infinitely many solutions (see
example 1 below).

Since the Blaschke product $\F_{n-2}/\F_{n-2}^*$ of order $n-2$
cannot take the same value on the $n$-point set $\S(\cC)$, there
always exists such a pair $\z_1\ne\z_2$ in $\S(\cC)$, that
(\ref{1.10}) holds, so $\a_{n-2}$ is uniquely determined.

\bigskip

The general MISP for CMV matrices we study here looks as follows.
Let $n=2l$ be even. Given first $n-m-1$ Verblunsky parameters
$\a_0,\ldots,\a_{n-m-2}$, and $2m$ eigenvalues
$\z_1,\ldots,\z_{2m}$, $1\leq m \leq n/2 = l$, find the rest $m$
parameters $\a_{n-m-1},\ldots,\a_{n-2}$ and thereby restore the
whole matrix $\cC$.
 \footnote{ $2m$ ``real'' parameters are given to find $m$
 ``complex'' ones.}
 Our main result
provides the conditions for this problem to have a unique solution.

Consider a pair of CMV matrices with the ``known'' parameters
$\cC(\a_0,\ldots,\a_{n-m-3};1)$ and
$\cC(\a_0,\ldots,\a_{n-m-2};1)$ and the systems of the \SZ\
polynomials
$$\{\F_0,\ldots,\F_{n-m-2};\wt\F_{n-m-1}\}, \qquad
\{\F_0,\ldots,\F_{n-m-1};\wt\F_{n-m}\},$$
 respectively. By the \SZ\ recurrences \eqref{1.0}
$$ \begin{aligned} \wt\F_{n-m-1}(z) &= z \F_{n-m-2}(z) -
\F_{n-m-2}^*(z),
\\
\F_{n-m-1}(z) &= z \F_{n-m-2}(z) - \bar\a_{n-m-2}\F_{n-m-2}^*(z),
\end{aligned}$$
so
\begin{equation} \label{1.11}
\F_{n-m-1}(z) - \wt\F_{n-m-1}(z) = (1-\bar\a_{n-m-2})\F_{n-m-2}^*.
\end{equation}
Similarly, for the pair $\cC(\l_0,\ldots,\l_{m-1};1)$ and
$\cC(\l_0,\ldots,\l_m;1)$ of ``unknown'' CMV matrices, $\l_j$ from
\eqref{1.71} with the Szeg\H{o} polynomials
$\{\L_0,\ldots,\L_m;\wt\L_{m+1}\}$ and \newline
$\{\L_0,\ldots,\L_{m+1};\wt\L_{m+2}\}$, respectively, one has
\begin{equation}\label{1.12}
\L_{m+1}(z)-\wt\L_{m+1}(z) = (1-\bar\l_m) \L_m^*.
\end{equation}
Now write (\ref{1.8}) with $k=m+1$:
\begin{equation}\label{1.13}
\frac{y_{m}(\z_j)}{y_{m+1}(\z_j)} =
\frac{x_{n-m-1}(\z_j)}{x_{n-m-2}(\z_j)}\,, \qquad j=1,\ldots,n,
\end{equation}
and observe that the right hand side of (\ref{1.13}) is known for
$j=1,2,\ldots,2m$. Indeed, let, e.g., $m$ be odd (for even $m$ the
calculation is the same). Then by (\ref{1.2})
$$x_{n-m-2}(z) = z^{-\frac{n-m-1}{2}} \k_{n-m-2}\F_{n-m-2}^*(z);
\quad x_{n-m-1}(z) = z^{-\frac{n-m-1}{2}} \k_{n-m-1}\F_{n-m-1}(z),
$$
so in view of (\ref{1.11})
$$\begin{aligned}
\frac{x_{n-m-1}(\z_j)}{x_{n-m-2}(\z_j)} &=
\frac{\k_{n-m-1}}{\k_{n-m-2}}
\frac{\F_{n-m-1}(\z_j)}{\F_{n-m-2}^*(\z_j)} \\
&= (1-|\a_{n-m-2}|^2)^{-1/2}
\frac{\wt\F_{n-m-1}(\z_j)+(1-\bar\a_{n-m-2})\F_{n-m-2}^*(\z_j)}{\F_{n-m-2}^*(\z_j)}\\
&= \rho_{n-m-2}^{-1} \left\{
\frac{\wt\F_{n-m-1}(\z_j)}{\F_{n-m-2}^*(\z_j)} + 1-\bar\a_{n-m-2},
\right\}, \quad \rho_i=(1-|\a_i|^2)^{1/2}.
\end{aligned}
$$
In the same way
$$ y_{m}=z^{-\frac{m+1}{2}}\k_{m,r}\L_m^*(z),
\quad y_{m+1}=z^{-\frac{m+1}{2}}\k_{m+1,r}\L_{m+1}(z),
$$
and with
$$ \k_{m,r}=\prod_{j=0}^{m-1} (1-|\l_j|^2)^{-1/2}=
\prod_{j=n-m-1}^{n-2} (1-|\a_j|^2)^{-1/2} $$ we have
$$ \begin{aligned}
\frac{y_{m}(\z_j)}{y_{m+1}(\z_j)} &= \frac{\k_{m,r}}{\k_{m+1,r}}
\cdot \frac{\L_{m}^*(\z_j)}{\L_{m+1}(\z_j)}=
 \rho_{m,r} \left\{
\frac{\wt\L_{m+1}(\z_j)+(1-\bar\l_m)\L_{m}^*(\z_j)}{\L_m^*(\z_j)}\right\}^{-1}\\
 &= \rho_{m,r} \left\{  \frac{\wt\L_{m+1}(\z_j)}{\L_m^*(\z_j)} +
1 + \a_{n-m-2} \right\}^{-1}, \quad \rho_{m,r} =(1-|\l_m|^2)^{-1/2}
=\rho_{n-m-2}.
\end{aligned}
$$
Using \eqref{1.13}, we end up with the following equalities for
$j=1,2,\ldots,2m$
\begin{equation}\label{1.14}
\frac{\wt\L_{m+1}(\z_j)}{\L_m^*(\z_j)} = -1 - \a_{n-m-2} + \frac{1
- |\a_{n-m-2}|^2}{ \frac{\wt\F_{n-m-1}(\z_j)}{\F_{n-m-2}^*(\z_j)}
+ 1 - \bar\a_{n-m-2}}\,.
\end{equation}

As the last step, we express the ratios in terms of the Weyl
functions (more precisely, their reciprocals)
$$ W(z)= \frac1{w(z)} = \frac{\wt\F_{n-m-1}(z)}{\F_{n-m-2}(z)}\,, \qquad
W_r(z)=\frac1{w_r(z)} = \frac{\wt\L_{m+1}(z)}{\L_m(z)}
$$
of the known $\cC(\a_0,\ldots,\a_{n-m-3};1)$ and unknown
$\cC(\l,\ldots,\l_{m-1};1)$, respectively. Indeed, for $|z|=1$
$$ \begin{aligned}
\wt\F_{n-m-1}(z) &= \prod_{i=1}^{n-m-1} (z-z_i), \quad |z_i|=1,
\quad
\wt\F_{n-m-1}(0)=(-1)^{n-m-1}\prod_{i=1}^{n-m} z_i=-\bar\b=-1, \\
\ov{\wt\F_{n-m-1}(z)}
&=\prod_{i=1}^{n-m-1}(z^{-1}-z_i^{-1})=-z^{-n+m+1}\wt\F_{n-m-1}(z),
\end{aligned} $$
so $\ov{\wt\F_{n-m-1}(\z_j)} = - \z_j^{-n+m+1} \wt\F_{n-m-1}(\z_j)$,
$$ \frac{\wt\F_{n-m-1}(\z_j)}{\F_{n-m-2}^*(\z_j)} =
- \frac{\z_j^{n-m-1} \overline{\wt\F_{n-m-1}(\z_j)}} {\z_j^{n-m-2}
\overline{\F_{n-m-2}(\z_j)}} = - \z_j\ov{W(\z_j)},
$$
and, similarly,
$$\frac{\wt\L_{m+1}(\z_j)}{\L_m^*(\z_j)}=-\z_j\ov{W_r(\z_j)}.
$$

Finally, we come to the following interpolation problem for the
Weyl function of the ``unknown'' CMV matrix
$\cC(\l_0,\ldots,\l_{m-1};1)$
\begin{equation}\label{1.15}
W_r(\z_j) = \z_j \frac{\ds(1+\bar{\a}_{n-m-2})\left\{ 1 - \a_{n-m-2}
- \bar\z_j W(\z_j) \right\} -
(1-|\a_{n-m-2}|^2)}{\ds\left\{1-\a_{n-m-2}-\bar\z_j W(\z_j)\right\}}
=:\o_j \,,
\end{equation}
$j=1,2,\ldots,2m$, or
\begin{equation}\label{1.16}
P^{(1)}(\z_j) - \o_j P^{(2)}(\z_j) =0, \quad j=1,\ldots,2m,
\end{equation}
with $\o_j$ defined in (\ref{1.15}), which we have denoted by
$(I_{2m})$ in the previous section. Now $\o_j\ne\infty$ since all
zeros of $\L_m$ are in the open unit disk $\dD$. The above
argument shows that \eqref{1.16} has a nontrivial solution
\footnote{Again, we assume that $\cC=\cC(\a_0,\ldots,\a_{n-2};1)$
with given $\a_0,\ldots,\a_{n-m-2}$ and the eigenvalues
$\z_1,\ldots,\z_{2m}$ does exist}
$$
\l= \begin{pmatrix} \L^{(1)} \\ \L^{(2)}\end{pmatrix} =
\begin{pmatrix} \wt\L_{m+1}\\\L_m\end{pmatrix}
$$
and
\begin{equation}\label{1.161}
h(\l)= 2m+2
\end{equation}

\begin{proposition} \label{propos1.1}
For the problem \eqref{1.16}\  $h(I_{2m})\ge 2m-1$.
\end{proposition}

{\it Proof.} Let $r$ be the minimal generator of \eqref{1.16}, and
suppose that $h(r)\le 2m-2$. By Theorem \ref{theor2.7} for the
second generator $q$ one has $h(q)\ge 2m+3$. It follows now from
Theorem \ref{theor2.8} and \eqref{1.161} that $\l=Sr$ with
 $\deg S\ge 2$, so
$$ \wt\L_{m+1}(z)=S(z)R^{(1)}(z), \qquad
\L_{m}(z)=S(z)R^{(2)}(z), $$ which is impossible, for $\wt\L_{m+1}$
and $\L_m$ have no common zeros. \hfill $\square$

There is some more information available about the solution $\l$.
Specifically,
\begin{equation}\label{1.17}
\deg \wt\L_{m+1}=m+1, \qquad \deg \L_{m}=m
\end{equation}
and
\begin{equation}\label{1.18}
 \wt\L_{m+1}(0)=-1.
\end{equation}

\medskip

In view of Theorem \ref{theor2.6}, it is easy to conclude from the
Proposition \ref{propos1.1} that if the data of interpolation
problem \eqref{1.16} correspond to a CMV matrix, then either
$h(I_{2m})=2m$ or $h(I_{2m})=2m-1$.

\begin{theorem}\label{theormain}
Let for the minimal generator $r$ of problem $(\ref{1.16})$
$R^{(1)}(0) \ne 0$ holds. Then $(\ref{1.16})$ has a unique solution
and, hence, the solution of the MISP is unique.
\end{theorem}

{\it Proof.} \ By Theorem \ref{theor2.8}
$$ \begin{aligned} \L^{(1)}(z) &= S(z)R^{(1)}(z) + T(z)
Q^{(1)}(z),\\
\L^{(2)}(z) &= S(z)R^{(2)}(z) + T(z) Q^{(2)}(z), \end{aligned} $$
where $r$ and $q$ are the minimal and second generators for the
problem (\ref{1.16}), respectively. Proposition \ref{propos1.1}
reads that either $h(r)=2m$ or $h(r)=2m-1$, so by Theorem
\ref{theor2.7} either $h(q)=2m+1$ or $h(q)=2m+2$. In the first case
$$ \deg R^{(1)} =m, \quad \deg R^{(2)}\leq m-1, \quad \deg
Q^{(1)}\le m, \quad \deg Q^{(2)}=m, $$ and in the second one
$$ \deg R^{(1)} \le m-1, \quad \deg R^{(2)}=m-1, \quad \deg
Q^{(1)}= m+1, \quad \deg Q^{(2)}\le m. $$
 In view of the degrees of $\L^{(1)}$ and $\L^{(2)}$, in both cases
$$ \deg S(z)=1, \qquad \deg T(z) =0,
$$
i.e., for $\L^{(j)}$ we have
\begin{equation} \label{1.18.5}
\begin{aligned}
\L^{(1)}(z) &= (az + b)R^{(1)}(z) + c
Q^{(1)}(z), \\
\L^{(2)}(z) &=  (az + b)R^{(2)}(z) + c Q^{(2)}(z), \end{aligned}
\end{equation}
For fixed $r$ and $q$, it is easy to see that $a$ and $c$ are
uniquely determined by condition (\ref{1.17}). If $R^{(1)}(0)\ne 0$,
then $b$ is also uniquely determined by (\ref{1.18}), but if
$Q^{(1)}(0)=0$, then problem has either no solutions, or infinitely
many solutions. But since the data of the problem is taken from the
CMV matrix, the solution does exist, so if $R^{(1)}(0)=0$, the
interpolation problem has infinitely many solutions and, hence, the
MISP may have infinitely many solutions. \hfill $\square$

\bigskip

In view of this theorem, two natural questions arise.

1) Is it possible for the minimal generator to have $R^{(1)}(0)=0$
when \eqref{1.15} is related to MISP for a certain CMV matrix?

2) Is it possible for the MISP to have more then one solution if
$R^{(1)}(0)=0$?

Both answers are positive, so in some special cases MISP with
non-unique solutions does exist, although the number of the
``pieces of information'' in the inverse data is equal to the
number of parameters to reconstruct. We provide examples for both
possible cases $h(I_{2m})=2m-1$ (example 1) and $h(I_{2m})=2m$
(example 2).

\medskip

{\bf Example 1.} \ Let $-1<b<1$ and $\cC=\cC(0,0,b;1)$. By the
Szeg\H{o} recurrences
$$ \F_1(z) = z, \ \F_2=z^2, \ \F_3=z^3-b, \ \F_3^*(z) = -bz^3 +1,
$$
and
$$ \wt\F_4 (z) = z \F_3(z) - \F_3^*(z) = (z^2-1)(z^2+bz+1),
$$
so the eigenvalues are
$$\z_{1,2} = \pm 1, \quad \z_{3,4} = - \frac{b}{2}
\pm i \sqrt{1-\frac{b^2}{4}}.
$$
We see that the pair $\z_1, \z_2$ does not determine $b$ uniquely,
although any other pair does.

However, the MISP of (non-unique) reconstruction of $b$ by the two
eigenvalues $\z_{1,2} = \pm 1$ is still possible. Find the
right-hand side of \eqref{1.15} for this case. First, consider the
Weyl function of the ``known'' left matrix $\cC(0;1)$. Its
Szeg\H{o} polynomials are $\F_1(z) = z$ and
$\wt\F_2(z)=z\F_1-1\cdot\F_1^*=z^2-1$. So, the reciprocal of its
Weyl function is $ W(z)={\ds \frac{z^2-1}{z}}$ and the right-hand
side of \eqref{1.15} is
$$
\o_j=\z_j\cdot
\frac{(1+0)\{1-0-\bar\z_jW(\z_j)\}-(1-0)}{1-0-\bar\z_jW(\z_j)}=
\z_j\cdot\frac{-\bar\z_jW(\z_j)}{1-\bar\z_jW(\z_j)}.
$$
For $\z_j=\pm1$ we have $\o_j =0$, so, according to \eqref{1.15},
the reciprocal of the Weyl function of the matrix
$\cC_r=\cC(-b;1)$ satisfies
$$ W_r(\pm1)=0.$$
Remind that it also must satisfy additional conditions
\eqref{1.17} and \eqref{1.18} with $m=2$.

So, following the procedure of solving the MISP, described above, to
reconstruct the inverse of the Weyl function of $\cC(-b;1)$, we need
to reconstruct a rational function ${\ds \frac{P^{(1)}}{P^{(2)}}}$
such that
\begin{equation}\label{1.30}
\begin{aligned}
\frac{P^{(1)}(z)}{P^{(2)}(z)}|_{z=\pm 1} =0, \quad
\deg P^{(1)}=2, \quad \deg P^{(2)}(z) =1,\\
P^{(1)}, P^{(2)} \ {\rm are \ monic \ and} \ P^{(1)}(0)=1.
\end{aligned}
\end{equation}

The corresponding interpolation problem for the vector-functions is
\begin{equation}\label{1.31}
P^{(1)}(\pm 1)=0.
\end{equation}
According to Proposition \ref{propos1.1}, the minimal generator of
the problem (\ref{1.31}) must have the height $\ge 1$, and according
to Theorem \ref{theor2.6}, it must have the height $\leq 2$. In
fact, it can be immediately checked that the non-trivial
vector-function of minimal height, corresponding to this problem, is
$\begin{pmatrix} 0\\ 1
\end{pmatrix}$, whose height is 1. Further, since the height of the minimal
generator is $1$ and there is $2$ points of interpolation, according
to Theorem \ref{theor2.7}, the second generator must have the height
$4$ and there is no solutions of height $2$. In fact, it is evident
that
the vector-polynomial $\begin{pmatrix} (z+1)(z-1)\\
0 \end{pmatrix}$, whose height is 4, solves (\ref{1.31}), and there
is no solutions of height $2$. Finally, the general solution of
(\ref{1.30}), is
$$\frac{P^{(1)}(z)}{P^{(2)}(z)}
=\frac{1\cdot (z+1)(z-1)+(z + a)\cdot 0}{1\cdot 0 + (z + b) \cdot 1}
= \frac{(z+1)(z-1)}{z + b},
$$
with arbitrary number $b$ (cf. \eqref{1.18.5}). However, only those
solutions with additional condition $|b|<1$ give us not only a
solution of (\ref{1.30}), but also the Weyl function of a CMV matrix
of the type $\cC(-b;1)$. In fact, let us find directly the Weyl
function of $\cC(-b;1)$. Its Szeg\H{o} polynomials are:
$$ \left\{\begin{array}{l}
\L_1(z)=z+b; \\
\wt\L_2(z)=z \L_1(z) - \L_1^*(z) = z(z+b) - (bz+1)= (z+1)(z-1).\\
\end{array}
\right.
$$
So,
$$W_r(z)=\frac{\wt\L_2}{\L_1}= \frac{(z+1)(z-1)}{z + b},
$$
as was to be checked.

\medskip

{\bf Example 2.} \  Consider a family of CMV matrices of order $4$:
$\cC(0,-y,-x;1)$; $-1<x,y<1$, and analyze the MISP of reconstruction
of the unknown $x, y$ by the four eigenvalues. Calculate for them
the Szeg\H{o} polynomials:
$$ \F_1(z)=z, \quad \F_1^*(z)=1;
$$
$$\F_2=z\F_1+y\F_1^*=z^2+y, \quad \F_2^*=1+z^2y;
$$
$$ \F_3=z\F_2+x\F_2^*=z(z^2+y)+x(1+z^2y)=z^3+xyz^2+yz+x, \quad
\F_3^*=1+xyz+yz^2+xz^3;
$$
$$\begin{aligned} \wt\Phi_4(z) &=
z\F_3-1\cdot\F_3^*=z(z^3+xyz^2+yz+x)-1-xyz-yz^2-xz^3
\\
&= (z^4-1) + (xy-x)z(z^2-1) = (z^2-1)(z^2+(xy-x)z+1).\end{aligned}
$$
Introducing the notation
\begin{equation} \label{1.19}
k:= xy-x,
\end{equation}
we express the eigenvalues of $\cC$ as
\begin{equation} \label{1.20}
\Sigma: \ \quad \z_{1,2}=\pm1; \quad \z_{3,4}=-\frac{k}{2}\pm i
\sqrt{1-\frac{k^2}{4}}; \qquad \z_3\ne \z_4 \ {\rm and} \ \z_4=\bar
\z_3.
\end{equation}

Hence, if $x$ and $y$ are related by (\ref{1.19}) and $k$ is fixed,
we have an infinite family of CMV matrices $\cC(x,y,0;1)$ with the
same spectrum $\Sigma$ (\ref{1.20}). According to the general theory
the auxiliary Weyl functions $W_r(z)$ and $W(z)$ of the matrices
$\cC_r=\cC(x,y;1)$ and $\cC(1)$, resp., take the same values on
$\Sigma$ for different $x$ and $y$, related by (\ref{1.19}). Find
the right-hand side of \eqref{1.15} for this case. First, consider
the Weyl function of the ``known'' left matrix $\cC(1)$. Its
Szeg\H{o} polynomials are $\F_0(z) = 1$ and
$\wt\F_1(z)=z\F_0-1\cdot\F_0^*=z-1$. So, the inverse of its Weyl
function is $ W(z)=z-1$ and the right-hand side of \eqref{1.15} is
$$
\o_j=\z_j\cdot
\frac{(1+0)\{1-0-\bar\z_jW(\z_j)\}-(1-0)}{1-0-\bar\z_jW(z_j)}=
\z_j\cdot\frac{-\bar\z_j (\z_j-1)}{1-\bar\z_j(\z_j-1)}=\z_j(1-\z_j).
$$
Let us directly check that the left-hand side of \eqref{1.15}
coincide with the obtained numbers. The Szeg\H{o} polynomials of
$\cC_r=\cC(x,y;1)$ are
$$ \L_1(z)=z-x, \quad \L_1^*(z)=1-xz;
$$
$$ \L_2(z)=z\L_1-y\L_1^*=z^2+(xy-x)z-y,
 \quad \L_2^*(z)=-yz^2+(xy-x)z+1;
$$
$$\wt\L_3(z)=z\L_2-\L_2^* = z^3 + (xy-x)z^2-yz+yz^2-(xy-x)z-1=
(z-1)(z^2+(k+y+1)z+1);
$$
$$ W_r(z) = \frac{\wt\L_3(z)}{\L_2(z)} =
\frac{(z-1)(z^2+(k+y+1)z+1)}{z^2+kz-y};
$$
$$W_r(1)=0 = \o_1;
$$
$$W_r(-1) = \frac{-2(2-(k+y+1))}{1-k-y}=-2 = \o_2;
$$
and, since $\z_j^2+k\z_j+1=0$; $j=3,4$, we have
$$ W_r(\z_j) = \frac{(\z_j-1)(y+1)\z_j}{-(y+1)} = \z_j(1-\z_j) = \o_j; \quad
j=3,4.
$$
If we solved the MISP following the procedure described above, we
have to find (non-uniquely) the inverse Weyl function $W_r(z)$ from
its value in the four eigenvalues and the additional conditions for
the numerator and denominator. In this example we will restrict
ourselves by illustrating that the minimal generator of the
corresponding interpolation problem does not satisfy the conditions
of theorem \ref{theormain}, which actually cause the existence of
infinitely many solutions.

Consider the interpolation problem, corresponding to this case:
$$
\left\{\begin{array}{l} R^{(1)}(\z_j) - \o_j R^{(2)}(\z_j); \quad
j=1,2,3,4;\quad \z_j\in\Sigma; \\
\o_1=0, \quad \o_2=-2, \quad \o_{3,4}=\z_{3,4}(1-\z_{3,4}).\\
\end{array}
\right.
$$
(since $\z_j\ne1$, we have $\o_j\ne 0$). We are looking for the
solution of height $4$: $R^{(1)}=z^2+\a z +\b$, $R^{(2)}=\g z + \d$,
so
$$
\left\{\begin{array}{l} 1+\a + \b =0; \\
1-\a+\b + 2(-\g+\d) =0;\\
\z_j^2+\a \z_j +\b -\o_j (\g \z_j +\d) =0, \quad j=3,4.
\end{array}
\right.
$$
Since $\z_j^2+k\z_j+1=0$, the last $2$ equations can be rewritten as
$$(\a-k)\z_j+\b-1-\o_j(\g \z_j + \d) = 0, \quad (\a-k)\z_j - \a -2 -
\o_j(\g \z_j +\d) =0; \ j=3,4.
$$
Exclude from these equations first $\g$, then $\d$:
$$
{\rm 1)} \ \left\{\begin{array}{l}
(\a-k)\z_3\o_4-(\a+2)\o_4-\g \z_3\o_3\o_4-\d\o_3\o_4 =0; \\
(\a-k)\z_4\o_3-(\a+2)\o_3-\g \z_4\o_4\o_3-\d\o_3\o_4 =0 .\\
\end{array}
\right.\Rightarrow
$$
$$(\a-k)(\z_3\o_4 - \z_4\o_3) - (\a+2)(\o_4-\o_3) - \g |\o_3|^2
(\z_3-\z_4) =0.
$$
$$(\a-k)(\z_3-\z_4) - (\a+2)(\o_4-\o_3) = \g |1-\\z_3|^2 (\z_3-\z_4);
$$
$$\g = \frac{\a-k}{|1-\z_3|^2} + (a+2)
\frac{\o_3-\o_4}{(\z_3-\z_4)|1-\z_3|^2} = \frac{\a-l}{|1-\z_3|^2} +
(\a+2)\frac{(1+k)}{|1-\z_3|^2} = \frac{(\a+1)(k+2)}{|1-\z_3|^2}.
$$
But $|1-\z_3|^2=2-\z_3-\bar \z_3 = 2+k \Rightarrow \g=\a+1$.
$$
{\rm 2)} \ \left\{\begin{array}{l}
(\a-k)\z_3\cdot \z_4\o_4-(\a+2)\z_4\o_4-\g \z_3\o_3 \z_4\o_4-\d\o_3 \cdot \z_4 \o_4 =0; \\
(\a-k)\z_4\cdot \z_3\o_3-(\a+2)\z_3\o_3-\g \z_4\o_3 \z_4\o_4-\d\o_4 \cdot \z_3 \o_4 =0 .\\
\end{array}
\right.\Rightarrow
$$
$$(\a-k)(\o_4 - \o_3) - (\a+2)(\z_4\o_4-\z_3\o_3) - \d |1-\z_3|^2
(\z_4-\z_3) =0,
$$
$$ \z_4\o_4 - \z_3\o_3 = \z_4^2(1-\z_4) -
\z_3^2(1-\z_3)=(\z_4-\z_3)(\z_4+\z_3-\z_4^2-\z_3^2-1) =
(\z_4-\z_3)(-k^2-k+1),
$$
$$ \o_4-\o_3 = (\z_4-\z_3) (1-\z_3-\z_4) = (\z_4-\z_3) (1 + k).$$
Hence,
$$ (\a-k)(1+k)+(\a+2)(k^2+k-1) = \d (2+k);
$$
$$ \a k (k+2) +(k-1)(k+2) = \d (k+2) \Rightarrow \d=\a k+k-1.
$$

3) We have from the second equation $\a=\d-\g \Rightarrow \a=\a k +
k -2 \Rightarrow$
$$ \a=-1; \quad \b=0; \quad \g=0 ; \quad \d=-1,
$$
$$ R^{(1)}(z)=z^2-z; \quad R^{(2)}(z)=-1,
$$
It can be immediately checked that the solution of height $3$ (such
that $R^{(1)} =\a z +\b$, $R^{(2)} = z+\g$ does not exist, so
$ r=\begin{pmatrix} R^{(1)}\\
R^{(2)}
\end{pmatrix} $ is the minimal solution. As we see, $R^{(1)}(0) =0$.

Finally, rewriting
$$ W_r(z) =
\frac{(z-1)(z^2+(k+1)z+1)+ (z^2-z)\cdot y}{z^2+kz+ (-1)\cdot y};
$$
we see in the non-determined (arbitrary) up to a constant factor
terms of the numerator $(z^2-z)\cdot y$ and denominator $-1\cdot y$
the components of the minimal generator (cf. \eqref{1.18.5} with
$b=y$).

\bigskip

{\it Remark.} \ Assume that $2m+1$ eigenvalues
$\z_1,\ldots,\z_{2m},\z_{1m+1}$ are known. Since
$$ w_2(z) = \frac{\Psi_m(z)}{\wt\Psi_{m+1}(z)} =
\frac{z^m+\ldots}{z^{m+1}+\ldots}
$$
the interpolation problem $w_2(\z_j)=\O_j$; $j=1,2,\ldots,2m+1$ has
obviously the unique solution.


\begin{thebibliography}{}


\bibitem{CMV2}
M. J. Cantero, L. Moral, and L. Vel\'azquez, Five-diagonal matrices
and zeros of orthogonal polynomials on the unit circle, Lin. Algebra
Appl. {\bf 362} (2003), 29--56.

\bibitem{CMV3}
M. J. Cantero, L. Moral, and L. Vel\'azquez, Minimal representations
of unitary operators and orthogonal polynomials on the unit circle.
Lin. Algebra Appl. {\bf 408} (2005), 40--65.

\bibitem{CMV4}
M. J. Cantero, L. Moral, and L. Vel\'azquez, Measures on the unit
circle and unitary truncations of unitary operators, J. Approx.
Theory {\bf 139} (2006), 430--468.

\bibitem{Chu-Golub} M.T.Chu, G.H.Golub. Structured eigenvalue problems.
Acta Numer., {\bf 11} (2002), 1--71.

\bibitem{GS} F. Gesztesy, and B. Simon, M-functions and inverse spectral
analysis for finite and semi-infinite Jacobi matrices, Journal
d'Analyse Math\'ematique, {\bf 73} (1997), 267--297.

\bibitem{GK2} L. Golinskii, and M. Kudryavtsev,  An inverse spectral
theory for finite CMV matrices, Preprint arXiv:0705.4353, 2007.

\bibitem{HL} H.Hochstadt and B.Lieberman, An inverse Sturm-Liouville
problem with mixed given data, SIAM J.App;.Math. {\bf 34}, pp.
676-680 (1978).

\bibitem{H} H.Hochstadt, On the construction of a Jacobi matrix from
mixed  given data, Linear Alg. Appl. {\bf 28}, pp. 113-115 (1979).

\bibitem{KN}
R. Killip, and I. Nenciu, CMV: the unitary analogue of Jacobi
matrices, Preprint arXiv:math.SG/0508113, 2005


\bibitem{Kud1} M.Kudryavtsev, The direct and inverse problem of
spectral analysis for five-diagonal symmetric matrices, I, Mat. fiz,
anal, geom (1998) {\bf vol. 5}, ¹3/4. pp 182--202.

\bibitem{Kud2} M.Kudryavtsev, The direct and inverse problem of
spectral analysis for five-diagonal symmetric matrices, II, Mat.
fiz, anal, geom (1998) {\bf vol. 6}, ¹1/2. pp. 55--80.

\bibitem{simA}
B. Simon, {\it Orthogonal Polynomials on the Unit Circle}, V.1:
Classical Theory, AMS Colloquium Series, American Mathematical
Society, Providence, RI, 2005.

\bibitem{simB}
B. Simon, {\it Orthogonal Polynomials on the Unit Circle}, V.2:
Spectral Theory, AMS Colloquium Series, American Mathematical
Society, Providence, RI, 2005.

\bibitem{sim2}
B. Simon, CMV matrices: five years after, to appear in Proceedings
of the W. D. Evans' 65th Birthday Conference.

\bibitem{sim3}
B. Simon, Rank one perturbations and the zeros of paraorthogonal
polynomials on the unit circle, to appear in J. Math. Anal. Appl.

\end{thebibliography}
\end{document}